\documentclass[
  a4paper, 
  reqno, 
  oneside, 
  11p
]{amsart}

\usepackage[utf8]{inputenc}

\usepackage[dvipsnames]{xcolor} % Must come before tikz!
\colorlet{cite}{LimeGreen!50!Green}
\usepackage{tikz}  % for diagrams
\usetikzlibrary{arrows,positioning}
\tikzset{ 
  baseline=-2.3pt,
  text height=1.5ex, text depth=0.25ex,
  >=stealth,
  node distance=2cm,
  mid/.style={fill=white,inner sep=2.5pt},
}

\usepackage{lmodern}
\usepackage{tikz-cd}
\usepackage{amsthm, amssymb, amsfonts}

\usepackage{graphicx,caption,subcaption}
\usepackage{braket}  % For nice sets
\usepackage{microtype}
\usepackage[%
  bookmarks=true,			% show bookmarks bar?
  unicode=true,			% non-Latin characters in Acrobat’s bookmarks
  pdftitle={Deformations of noncompact manifolds},		%
  pdfauthor={Elizabeth Gasparim}{},	% 
  pdfkeywords={Moduli space}{vector bundle}{}{},	% list of keywords
  colorlinks=true,		% false: boxed links; true: colored links
  linkcolor=Blue,			% color of internal links
  citecolor=cite,		% color of links to bibliography
  filecolor=magenta,		% color of file links
  urlcolor=RoyalBlue			% color of external links
]{hyperref}				% One of the most useful packages I have found.
\usepackage{cleveref}

\newtheoremstyle{mydef}
  {}		% Space above environment
  {}		% Space below environment
  {}		% Body font
  {}		% Indent amount (empty = no indent, \parindent = para indent)
  {\scshape}	% theorem head font
  {. }		% Punctuation after heading
  { }		% Space after heading
  {\thmname{#1}\thmnumber{ #2}\thmnote{ #3}}	% Heading spec

\theoremstyle{plain}	% 'plain' is the default.  The others are 'definition' and 'remark'.
\newtheorem{theorem}{Theorem} % putting [section] on the end here tells latex to number the theorem environment within sections, ie. Theorem 2.3 for the third theorem in section 2.
\newtheorem*{theorem*}{Theorem}

\theoremstyle{mydef} % Here we have used the custom theoremstyle defined above instead of the usual 'definition' style.

\newtheorem*{conjecture*}{Conjecture}
\theoremstyle{remark}

\newtheorem*{remark*}{Remark}

\newtheorem{example}[theorem]{Example}

\newtheorem*{proposition*}{Proposition}
\newtheorem*{lemma*}{Lemma}
\newtheorem*{corollary*}{Corollary}
\theoremstyle{definition}
\theoremstyle{remark}

\DeclareMathOperator{\Tot}{Tot}

\newcommand{\ce}{\mathrel{\mathop:}=}
\usepackage{abstract}
\author{Elizabeth  Gasparim}
\address{E. Gasparim - Depto. Matem\'aticas, Univ. Cat\'olica del Norte, Antofagasta, Chile. 
etgasparim@gmail.com}
\title[Ballico's moduli spaces and  applications ]{Some contributions of Edoardo Ballico to \\ Moduli spaces and their applications }

\begin{document}

\maketitle
\begin{abstract}
\it \phantom{xxxxxx}To Edoardo Ballico, the most interesting person I have ever met. \\ \\ \centerline{\large \bf Happy 70th Birthday!}
\end{abstract}

%\tableofcontents

\section{Introduction}\footnote{Section 9 was written by Wojciech Kucharz}
When I was offered the opportunity to write  about moduli spaces of vector bundles in the work of Ballico, 
my first doubt was how I would  transmit to a reader just how much fun it is to collaborate with him.
Luckily, the editors  Alessandra Bernardi and Claudio Fontanari, allowed us to start with an anecdote. 
Faced with the difficulty of deciding which anecdote to tell, since there are so many interesting stories about
Ballico, I ended up with 
two  representatives: the one of the {\it ridiculous} and of the one about {\it balancing}. These
appear at the start  of sections \ref{ridiculous} and \ref{balancing}.

\section{A Landau--Ginzburg model without projective mirrors}\label{ridiculous}
 Here is how we arrived at  the ridiculous. 
 Together with  collaborators, I had spent an entire semester trying to write 
an example for the  Homological Mirror Symmetry (HMS) Conjecture by Kontsevich.
 We were concentrated in the  
 case of the semisimple adjoint orbits of $\mathfrak{sl}(2, \mathbb C)$, viewed as a
symplectic Lefschetz fibration $\mathrm{LG}(2)$ via the standard potential function obtained by
 Cartan--Killing pairing with the fixed regular element $\mathrm{Diag}(1,-1)$.
  
HMS predicts that for every  (subvariety of) a projective variety $X$  there exists a 
(symplectic) Landau--Ginzburg model $(Y,f)$ 
satisfying the categorical equivalence 
\begin{equation} Fuk(Y,f) \equiv D^b(Coh \, X)\text{,}\end{equation}
where $Fuk(Y,f) $ denoted the Fukaya category of  vanishing cycles of $(Y,f)$.
We had calculated $Fuk(Y,f)$  for $Y=\mathrm{LG(2)}$ and our goal was to find the appropriate $X$. 
Since we participated in several HMS conferences,
we had asked a few experts, without getting any hint of an answer, except for the fact that many
people mentioned that it must be all very trivial and easy, because it is an example
coming from Lie theory. Yet, we could not find $X$.

 Having just arrived to visit Edoardo in Trento, I made my way to his office together with 
 two of my graduate students. Immediately and typically, as soon as we arrived he asked:\\
 
 Edo: -- What work are you doing now?\\
 
 Eli: -- I have a Fukaya--Seidel category, and I want you to give me an algebraic variety $X$ so that
 its derived category of coherent sheaves is isomorphic to this one.\\

 Then I wrote down this lemma on his blackboard:

\begin{lemma*}
\label{teosl2}   {$Fuk(\mathrm{LG}(2))$}
is generated by two Lagrangians $L_0$ and $L_1$ with morphisms: 
\begin{align*}
\label{morphismsfukaya}
\mbox{Hom} (L_i, L_j) 
\simeq
\begin{cases}
\mathbb Z \oplus \mathbb Z [-1] & i < j \\
\mathbb Z                       & i = j \\
0                               & i > j
\end{cases}
\end{align*}
where we think of $\mathbb Z$ as a complex concentrated in degree $0$ and $\mathbb Z [-1]$ as its shift, concentrated in degree $1$, and the products $m_k$ all vanish except for $m_2(\,\cdot\,, \mbox{id})$ and $m_2({\mbox{id},}\,\cdot\,)$.
\end{lemma*}

 He read it quickly, had a big laugh, and replied: \\
 
 Edo: -- Ridiculous! \\
 
 The students were aghast. I replied.\\
 
 Eli: -- Great! So, now we have a theorem! \\
 
 And we did have a theorem, his statement that it was a ridiculous question, which greatly amused him,
  simply meant that he could immediately see that such an $X$ did not exist. 
  
  \begin{remark*}Of course, 
  having instant clarity  about such a broad question is the result of an extensive knowledge 
  and experience about vector bundles, sheaves and their moduli on varieties of any type, for instance:
  rational surfaces \cite{87}, surfaces with negative Kodaira dimension \cite{92},  
  K3 surfaces \cite{93},  curves \cite{97a, 97b},  real algebraic surfaces \cite{00},
 threefolds containing ruled surfaces \cite{05}, Calabi--Yau threefolds \cite{21},   surfaces of general type \cite{22a}.
  \end{remark*}

  Edoardo's observation that the mirror we were looking for did not exist, allowed us  to 
  verify that we had the first 
  (and unique up to now) known example  
  of a symplectic LG model whose Fukaya category is not $DCoh(X)$ for any 
  algebraic variety $X$. Hence, we had constructed a {\it vampire} for 
  the HMS conjecture, by showing:

\begin{theorem*}\cite{18b} HMS fails for $\mathrm{LG}(2)$:
%[Thms.\thinspace 4.1,7.6]
\begin{itemize}
\item $\mathrm{LG}(2)$ has no projective mirrors.
\item $\overline{\mathrm{LG}(2)}$ has no projective mirrors.
\end{itemize}
That is,
for any (subvariety of a) projective variety $X$  we have 
 $${Fuk(\mathrm{LG}(2))} \not\equiv D^bCoh(X).   $$

\end{theorem*}

Once Edoardo had entered in the world of mirror symmetry, it was easy to get him
to show us how to solve other  question motivated by HMS. 

\section{Vector bundles, Lagrangian skeleta  and  duality}\label{dual}

We have several publications which arose from Edoardo's observations that
either certain varieties or maps are rational. 
In this subsection I mention two  results obtained using rational maps. 
The first one is a description of duality which interchanges birational maps with
deformations of complex structures. \\

During a  visit to Trento, Edoardo saw that I was calculating the {\it Lagrangian skeleton} 
of $T^*\mathbb P^n$.  By a Lagrangian skeleton, we mean the following.
 Let $L$ be
the union of the stable manifolds  of the Hamiltonian flow  with respect to a
K\"ahler metric; when $L$ is of middle dimension, it is called the 
Lagrangian skeleton.
We had calculated the entire skeleton for the case of $n=3$ and 
I was staring at it in a state of confusion, because I could not 
see any relation between the components of the skeleton. 
Here is the example:

\begin{example} The Lagrangian skeleton of $T^*\mathbb P^3$ consists 
of the following 4 submanifolds:
$$L_0= T^*_{[1,0,0,0]}\mathbb P^3\sim \mathbb C^3, \,\,  L_1=
  \mathcal O_{\mathbb P^1}(-1)\oplus \mathcal O_{\mathbb P^1}(-1), \,\, L_2=\mathcal O_{\mathbb P^2}(-1),
\,\,  L_3=  \mathbb P^3 .$$
\end{example}

Having taken a look at the example, Edoardo immediately noticed
the existence of birational transformations of $T^*\mathbb P^3$ taking $L_i$  to $L_{i+1}$. 
These gave rise to the  first vertical arrow in the commutative diagram appearing 
in the statement that follows.

  \begin{theorem*}\cite{23b} \label{t1}
  %[Thm.\thinspace 1.1]
    The following diagram commutes:
\begin{figure}[h]$$\label{diagr}
\begin{tikzcd}[row sep=10ex, column sep = -1.5ex,
/tikz/column 5/.append style={anchor=base east},
/tikz/column 6/.append style={anchor=base west}
]
  L_j \arrow[d,swap, "\textup{\bf bir}"] 
& \subset T^*\mathbb{P}^{n-1} \arrow[rrr, Leftrightarrow, "\textup{\textbf{dual}}"]   
&\phantom{xxxxx} & \phantom{xxxxx} & \phantom{.}
 \phantom{x...}\mathcal{O}_{Z_{n}^{\circ}}(j)
& \oplus 
&\!\! \mathcal{O}_{Z_{n}^{\circ}}(-j)\phantom{xxx} \\
 L_{j+1}
& \subset T^*\mathbb{P}^{n-1} \arrow[rrr,Leftrightarrow,swap, "\textup{\textbf{dual}}"]
& & & \phantom{.}
 %|[blue, rotate=-15]| 
\!\!\!\mathcal{O}_{Z_{n}^{\circ}}(j+1)
& \oplus \arrow[u, swap, "\textup{\bf def}"] 
& \mathcal{O}_{Z_{n}^{\circ}}(-j-1) .
\end{tikzcd}
$$
\caption*{Duality between Lagrangians and vector bundles.}
\end{figure}
\end{theorem*}

 A remarkable feature of this diagram and the duality it represents is the fact that 
 deformations of complex structures occur as dual to birational transformations. 
 Here, the second vertical arrow is deformation of rank 2 vector bundles on the 
 open surface $Z_n^o$. In this notation $Z_n \simeq \Tot(\mathcal O_{\mathbb P^1}(-n))$
 is the total space the line bundle $\mathcal O(-n)$ over $\mathbb P^1$ and $Z_n^o$ is its {\it collar}, 
 that is, the variety obtained from $Z_n$ after removing the zero section of $\mathcal O(-n)$. 
 In particular $Z_1^o \simeq \mathbb C^2 \setminus \{0\}$.
 We will revisit vector bundles on $Z_n$ in section \ref{balancing}. \\

Another case when Edoardo's constructions of rational maps were essential
 was when he formalised a way to look at some of our symplectic Landau--Ginzburg models 
 from  the algebro-geometric viewpoint \cite{17}. These LG models had been constructed
 completely with Lie theoretical methods, but we needed them to be compactified, by which 
 we mean not only compactify the variety, but also the map, which at first instance 
 can only be made rational. This happens because the LG models starts with a holomorphic map 
 $f\colon Y \rightarrow \mathbb C$, and the compactified version should give an extension of $f$ 
  to a map to  $\overline Y \dashrightarrow \mathbb P^1$, thus, necessarily only rational since $f$ is not constant. 

In the case when the LG model consists of 
$Y= \mathcal O_n$, which  is the minimal semisimple orbit of $\mathfrak{sl}(n+1,\mathbb C)$, 
and $f=f_H$ is the potential obtained by paring with $H$ via the Cartan--Killing form, the 
compactification may be described by the following diagram:
\begin{figure}[h]
	\centering
	\includegraphics[width=0.5\linewidth]{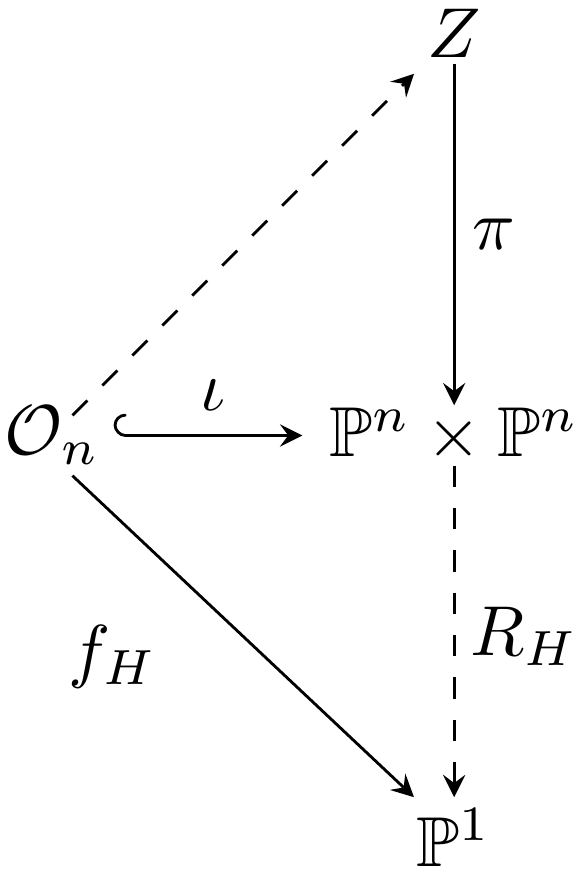}
	%\caption{}
	\label{fig:commutative-diagram-tikz-cd}
\end{figure}

\noindent where $Z$ is the closure of the graph of the blow-up map constructed 
for resolving the indeterminacy locus of $f_H$. Having obtained the desired compactification, 
we were then able to handle the question of Hodge theoretical invariants. 
 
 \section{Hodge theoretic invariants of Landau--Ginzburg models}

The original version of the Geometric Mirror Symmetry conjecture includes the statement  that for each Calabi--Yau 
variety $X$ there exists  a dual Calabi--Yau variety $X^\vee$ such that if the diamond of $X$ is given below, 
then the diamond of $X^\vee$ is obtained from the one of $X$ by reflection on the 45 degree line, 
that is, by exchanging  black and blue in the figure below, by flipping over the diagonal line passing thought the symbol for 
Serre duality. 

{\small \begin{displaymath}
\begin{array}{cccccccccc}
 & & & & h^{n,n} & & & & & \\
 & & & & & & & & \\
 & & & h^{n,n-1} & & h^{n-1,n} & & & & \\
 & & & & & & & & & \\
 & & h^{n,n-2} & & h^{n-1,n-1} & & \textcolor{Cyan}{h^{n-2,n}} & & & \\
 & & & &  & & & \textcolor{Cyan}{\ddots} & & \\
 h^{n,0} & & & {\cdots} & {\stackrel{\curvearrowleft}{_{Serre}}}{} &\textcolor{Cyan} {\cdots} & & & \textcolor{Cyan}{ h^{0,n}} & 
 \updownarrow_{Hodge \star}\\
 & \ddots & & & {} & & & & & \\
 & & h^{2,0} & & \textcolor{Cyan}{h^{1,1}} & & \textcolor{Cyan}{h^{0,2}} & & & \\
 & & & & & & & & & \\
 & & & \textcolor{Cyan}{h^{1,0}} & & \textcolor{Cyan}{h^{0,1}} & & & & \\
 & & & & & & & & & \\
 & & & & \textcolor{Cyan}{h^{0,0}} & & & & &\\
 & & & & \stackrel{\longleftrightarrow}{_{Conjugation}} & & & & &
\end{array}
\end{displaymath}
}
\vspace{3mm}

If a mirror $X^\vee$ failed to exist for $X$, then the variety $X$ would be said to behave like a {\it vampire}.
The validity of the conjecture would then imply that vampires do not exist. 

In the Homological version of the Mirror Symmetry Conjecture, there 
does exist an underlying statement about 
duality of Hodge theoretical invariants. Nevertheless, it was not enough to consider 
usual Hodge numbers of a variety. 

At this point we find that it is about time to 
clarify that by a Landau--Ginzburg model we mean a pair $(Y,f)$ of a variety
together with a complex function $f\colon Y \rightarrow \mathbb C $ (or $\mathbb P^1$).
If the  LG model appears on the $A$-side, $Y$ will be a symplectic variety, 
while on the $B$-side $Y$ will be an algebraic variety, but in either case it is 
typically noncompact. So, by applying the usual definitions of Hodge invariants, 
one would immediately obtain diamonds with infinities on it. 
Thus, one wishes to compactify $Y$. However, as we showed in 
\cite{18a}, the resulting Hodge diamond varies wildly according to the choice 
of compactification.

In 2017 Katzarkov, Kontsevich and Pantev defined 3 new Hodge theoretic 
invariants for LG model, which  considered not just the variety $Y$ but also 
the potential function $f$. They then conjecture the equality of these three invariants.
However, in 2018 Lunts and Przyjalkowski proved that the conjecture 
fails in complex dimension 2. Nevertheless, it is still an interesting 
question to verify which values are assumed by such invariants. We proved
their equality for the minimal adjoint orbit $\mathcal O_n$ endowed with the compactification
described in section \ref{dual}.

\begin{theorem*}\cite{23d}
%[Thm.\thinspace 6.1]
The $\mathrm LG$ model $(\mathcal O_n,f_H)$ admits a weakly tame compactification and satisfies the KKP conjecture.
\end{theorem*}

Whether the KKP conjecture holds true for other adjoint orbits remains an open question at this time.

\section{Existence of instantons and elementary transformations}\label{balancing}

When asking a question to Edoardo,  I got one of these three  reactions: 
\begin{itemize}
\item either he was 
not interested, and then the answer typically was \\

--Edo: Not my business. \\

\item  or  he liked the problem and would immediately  start working on it, in which case
new results would appear very fast,

\item or else  (this I saw just a small number times)  we got to the situation where he liked the problem, but did not 
have any immediate answer. In such a case, something amusing would happen. 
\end{itemize}

In 2002 we were discussing the problem of the existence of instantons, that is,  connections minimising the Yang--Mills 
functional
$$YM(A) = \int |F_A|^2 \textnormal{dvol}, $$
where $A$ is a connection and $F_A$ is the curvature of $A$.
 By the Kobayashi--Hitchin correspondence, instantons can be equivalently described by stable vector bundles,
whose  second Chern class corresponds to the instanton charge, i.e. there is a  1--1 correspondence
\[
\left\{
\begin{array}{c}
\text{irreducible $\mathrm{SU} (2)$-instantons} \\ 
\text{of charge $n$}
\end{array}
\right\}
\leftrightarrow
\left\{
\begin{array}{c}
\text{stable $\operatorname{SL} (2, \mathbb C)$-bundles} \\
\text{with $c_2 = n$}
\end{array}
\right\}.
\]\\
The open question of interest at that point was that of existence of instantons with prescribed {\it local charge} around a 
$-1$-line, that is an $\ell \simeq \mathbb P^1$ inside a surface having $\ell^2=-1$.
Let me first explain what is meant by  the  local charge. It is, in fact, a local version of the Euler characteristic, defined as follows.

If $\pi\colon (\widetilde{X}, \ell) \rightarrow (X,x)$ is a resolution of an isolated singularity, 
and $\mathcal F$ is a coherent sheaf of rank $n$ on $\widetilde{X}$, then the {\it local holomorphic Euler 
characteristic} of $\mathcal F$ is 
$$ \chi(\ell, \mathcal F) = h^0(x,Q) + \sum_{i=1}^n (-1)^ih^0(X, R^i\pi_*(\mathcal F)) $$
where $h^0$ is the dimension of the $0$-th \v{C}ech cohomology,
 $R^i\pi_*(\mathcal F)$ is the $i$-th higher derived image of $\mathcal F$  and $Q$ is 
the skyscraper sheaf supported at $x$ defined by the exact sequence
$$0 \rightarrow \mathcal \pi_*F \rightarrow (\pi_*\mathcal F)^{\vee\vee} \rightarrow Q \rightarrow 0.$$

% If $X$ is  a compact orbifold and $\mathcal F$ a sheaf over $X$, there exists 
%the global  {\it holomorphic orbifold Euler characteristic} of $\mathcal F$ 
%$$\chi_{orb}(X, \mathcal F) = \int_X  \mbox{ch}(\mathcal F) \mbox{td}(X).$$
%It satisfies:
%$$\chi(X, \mathcal F)= \chi_{orb}(X, \mathcal F) + \sum_{x \in \mbox{Sing}(X)}\mu(x,\mathcal F).$$ 
Therefore, for  $SU(2)$ instantons defined onver a 2-dimensional neighbourhood of  a $-1$-line, we needed to calculate:
$$ \chi(\ell, \mathcal F)= h^0(x,Q) + h^0(X, R^1\pi_*(\mathcal F)) $$ for rank 2  bundles with vanishing fist Chern class. For such computations, it 
was sufficient to consider  $\mathcal F$ as a rank 2 bundle over the blow up of $\mathbb C^2$ at the origin; in which 
case there exists an integer $j$, called the splitting type of $\mathcal F$, such that 
$\mathcal F\vert_{\ell} = \mathcal O_{\ell}(j) \oplus \mathcal O_{\ell}(-j)$. 
Having determined this $j$, I had obtained sharp bounds:
$$j \leq \chi(\ell, \mathcal F) \leq j^2,$$
but the question of whether intermediate values of $\chi$ occurred was entirely open, and it seemed out of reach, 
because,  we essentially needed to {\it guess} the appropriate vector bundle for each value of $\chi$.

But, Edoardo had become nervous at the sight of the interesting yet out of reach question, so we went for a walk 
at a nearby garden. As it turned out, the garden had a swing, and Edoardo sat on the swing and started moving back and forth, 
but he was a bit unbalanced, so it was not working very well. A couple of times, he stopped and started again, until
he got it to oscillate nicely. Then he was very amused, had a big laugh and said:\\

--Edo: Ah, now I am balanced! So,  I know. \\

Then he left. About a half an hour later, he sent me the proof of existence of bundles realizing every intermediate value of $\chi$, 
which he wrote as a game strategy of balancing  \cite{03}, starting with an unbalanced bundle, and  winning the game when the bundle 
had become balanced. 
It was a lovely way to solve the problem, inspired by playing at the garden.\\
 
 Some time later, together with T. K\"oppe \cite{09a, 09b} we studied the cases of local moduli of vector bundles on 
 a neighbourhood of  a $-k$-line in a surface, that is,  when $\ell\simeq \mathbb P^1$ with $\ell^2=-k$.
We determined sharp bounds for the local Euler characteristics and described
 the dimension and structure of these local moduli for each splitting type. 
 One interesting feature of the cases $k>1$ is that the lower bound for the local Euler characteristic is $k-1$, 
 from which it follows that when $k>2$ a neighbourhood of the $-k$ line does not support any instanton of charge $ 1$. 
 Hence the moduli of charge $1$  instantons on $\mbox{Tot}(\mathcal O_{\mathbb P^1}(-k))$ is empty whenever $k>2$. 
 Existence of bundles realising intermediate
 numerically admissible values of local charge however, remains an open problem for the cases of $-k$-lines when $k>2$.\\
 
 Still considering  questions about moduli spaces of vector bundles motivated by 
  instantons in mathematical physics, 
 together with C. Eyral  \cite{12}  we calculated
 the fundamental group of  moduli  of anti-self-dual connections on 
 blow-ups of $\mathbb P^2$, 
 proving that  the  fundamental group of such a moduli space is either $\mathbb Z$ of else $\mathbb Z_2$. 
 
 Edoardo, Rubilar and myself have also compiled a list of 25 open questions about vector bundles 
 and their moduli in \cite{22b}.

\section{Moduli of bundles on noncommutative Calabi--Yau varieties}

A subject that is always interesting to Edoardo is  deformation theory. In  \cite{23c} we expanded 
our search for deformations 
to include noncommutative ones, motivated by the notion of deformation quantization in mathematical physics. 
We then looked at the question of how moduli spaces of vector bundles on a Calabi--Yau threefold react to noncommutative 
deformations of the Calabi--Yau. We studied the case of the noncompact Calabi--Yau threefolds
$W_k = \Tot (\mathcal O_{\mathbb P^1}(-k) \oplus \mathcal O_{\mathbb P^1}(k-2))$. 
To calculate noncommutative deformations, we first described all Poisson structures on $W_k$ in \cite{23a},
representing the space of Poisson structure as a module over global functions. For each case, we named the 
generators of such modules as {\it basic} Poisson structures; these are somehow the 
simplest ones. It turns out that noncommutative 
deformations in the direction of basic Poisson structures are trivial in the sense that
the corresponding quantum moduli spaces of vector bundles are isomorphic to their classical limits 
on the commutative variety. 

On the other side of the Poisson spectrum live those Poisson structures that we named {\it extremal}.
A Poisson structure on $W_k$  is called extremal when it vanishes on the first formal neighbourhood of $\mathbb P^1$
inside $W_k$. 
Thomas K\"oppe had constructed  the  moduli spaces $\mathfrak M_j (W_k)$ of rank 2 bundles on $W_k$ with splitting type $j$
in his PhD thesis, and we described the corresponding moduli spaces of rank 2 bundles on  noncommutative deformations 
$\mathcal W_k$ of $W_k$. For the case when the first order term of the $\star$-product is determined by the 
Poisson structure $\sigma$,  we denote the moduli of rank 2 bundles on $\mathcal W_k$ with splitting type $j$ 
by $\mathfrak M_j^{\hbar} (\mathcal W_k, \sigma)$. We proved:

\begin{theorem*}\cite[Thm.\thinspace 5.4, 6.4]{23c} Let $k=1$ or $2$.
If $\sigma$ is an extremal Poisson structure on $W_k$, then
 the quantum moduli space $\mathfrak M_j^{\hbar} (\mathcal W_k, \sigma)$ can 
be viewed as the étale space of a constructible sheaf  $\mathcal E_k$ of generic 
rank $2j-k-1$ over the classical moduli space 
$\mathfrak M_j (W_k)$ with singular stalks of  all  ranks up to $4j -k-4$.
\end{theorem*}

It is entirely possible that this type of behaviour occurs in greater generality, 
including higher dimensions, but so far the general question remains open. 

\section{The number of components of moduli schemes of sheaves}

Let us now take a large jump backwards in time.
The three earliest writings of Edoardo about  moduli spaces of vector bundles on  surfaces consider those  surfaces which are
rational, then those with negative Kodaira dimension,
 and subsequently K3 surfaces.
His first result of this series is as follows.

\begin{theorem*}\cite{87}  If $S$ is a rational surface, then   there exists an ample divisor $H$ on $S$ such that for every value of $c_1$ and  $c_2$,
 the moduli scheme, $\mathfrak M_H(r,c_1,c_2)$ of $H$-stable vector bundles of rank $r$ with first Chern class $c_1$
  and second Chern class $c_2$ is smooth, irreducible, and unirational (if not empty). 
\end{theorem*}
  
For  surfaces with negative Kodaira dimension, Edoardo studied the birational structure of the exceptional irreducible components 
of  moduli schemes of  vector bundles \cite{98b}. In particular,  in 1992, he considered  the case of 
a ruled surface  over an algebraically closed field $k$ of characteristic $\neq 2$ with relatively minimal model $X\rightarrow S$ 
 obtained by blowing-up of  a geometrically ruled surface $S$, so that there are $\tau =K_S- K_X$ exceptional curves $D_i\subset X$. 
 Denoting by  $D_{\tau+2}$ be the pullback on $X$ of a ruling of $S$, the curves $D_i, i\leq \tau +2$ form a $\mathbb Z$-basis of $NS(X)$. 
 Let $\det=\sum d_iD_i$ [resp. $H=\sum h_iD_i$]  be a class [resp. an ample class] in $NS(X)$ (with $h_i,d_i\in \mathbb Z)$ and 
 denote by $\mathfrak M_X(2,(h_i),(d_i),c_2)$ the moduli space of rank 2 $H$-stable bundles on $X$ of determinant $\det$ and  second Chern class $c_2$.
 
\begin{theorem*}\cite{92} 
If $X$ and $X'$ are two such surfaces (associated with  blow-ups of the same surface $S$), then the moduli spaces
$\mathfrak M_X(2,(h_i),(d_i),c_2)$ and $\mathfrak M_{X'}(2,(h_i),(d_i),c_2)$ are birationally equivalent (and generically reduced of expected dimension) provided $c_2$ is large enough. 
 \end{theorem*}
 
 The following year,  Edoardo and Chiantini considered moduli  of stable bundles on $K3$ surfaces.
  Let $S$ be a minimal $K3$ surface, $H$ and  $L$ line bundles on $S$, with $H$ very ample satisfying $H\cdot L \leq 0$. 
  Let $\mathfrak M(L,c_2)$ denote the moduli space of rank 2 $H$-stable vector bundles on $S$ with determinant isomorphic to $L$ and second Chern class $c_2$.

 \begin{theorem*}\cite{93}   
    Let $E$ be an algebraic rank 2 vector bundle on $S$ with determinant isomorphic to $L$ and second Chern class $c_2$. 
If there exists an integer $m$ such that $\chi(E(mH))=1$ and  $mH+L$ is very ample, then there exists a birational map from 
$\mathfrak M(L,c_2)$ to the symmetric $(c_2+m^2H^2+mL\cdot H)$-product of $S$. It follows that in this case $\mathfrak M(L,c_2)$ is irreducible.
  \end{theorem*}  
  
  One year later, Edoardo obtained the following very strong result about the existence moduli 
  of spaces with infinitely many irreducible components.
  Let $k$ be an algebraically closed field. 
  Let $W$ be a smooth complete surface over $k$. 
  
   \begin{theorem*}\cite{94}  
  There exists a smooth surface $S$ obtained as the blowing-up of $W$ at finitely many points and 
  such that for all integers $c_2 \geq 2$ the moduli space of simple rank 2 vector bundles on $S$ 
  with determinant $\mathcal O_S$ and second Chern class $c_2$ has infinitely many components.
   \end{theorem*}  
   
   In contrast to this 1994  result  for surfaces, for the case of curves, Edoardo 
   found an explicit upper bound to the number of irreducible components of the moduli scheme of stable rank $r$ torsion-free sheaves 
   of fixed degree $d$ on an integral curve $X$. 
   Indeed, considering the moduli scheme of sheaves with fixed formal isomorphism type around singularities of the curve $X$, Edoardo
found   such  an upper bound to the number of components of the moduli scheme of torsion free sheaves, which is a linear function of $r$
     with coefficients depending on the singularities of the curve $X$, on  its arithmetic genus and on the number of components for rank-1 case. In particular, 
     he obtained the following result:
    
       \begin{theorem*}\cite{97a}   
The number of irreducible components of the moduli scheme of stable rank $r$ torsion-free sheaves 
   of fixed degree $d$ on an integral curve $X$ is finite and does not depend on $d$.
     \end{theorem*}

For the case of dimension 3, note also  that a lower bound for the number of components of the moduli schemes of stable rank  2  vector bundles on projective  threefolds
was given by Edoardo and Mir\'o-Roig in \cite{99}.

\section{Existence of singularities on moduli of vector bundles}

For simplicity, let me state here a result of Edoardo's about the existence of singularities in moduli of bundles on threefolds, 
and observe that, in \cite{98a}  he also obtained  similar results about the existence of singularities in moduli of bundles  over  $n$-folds, for any $n\geq 4 $.

  \begin{theorem*}\cite{98a}   
Assume characteristic 0. Fix a smooth projective threefold $X$, an $H\in \mathrm{Pic}(X)$ with $H$ ample and $d \in \mathbb Q$. 
Let $\mathfrak M(X, 2,\mathcal O, d, H)$ be the moduli scheme of rank 2 vector bundles, $E$, on $X$ with $\det E = \mathcal O_X$
and $c_2(E)\cdot H = dH^3$. Then for all integers $k\geq 2$ we have 
$$Sing\left( \mathfrak M(X, 2,\mathcal O, d, H)\right) \neq \emptyset .$$
     \end{theorem*}

\section{Moduli flexibility of   real  manifolds (written by Kucharz)}

I am grateful to  Wojciech Kucharz for providing  the text for this section. 

In \cite{13,14} Ballico and Ghiloni introduce and studied the notion of deformation of real algebraic manifolds.
Of course, the motivation came from the classical theory of deformation of complex analytic manifolds
and complex algebraic manifolds.
Roughly speaking, the main result of these papers can be summarized as the following principle.\\

{\sc Principle of real moduli flexibility}: The algebraic structure of every real algebraic manifold of positive 
dimension can be deformed by an arbitrary large number of effective parameters.\\

The possibility of an {\it arbitrary large number of effective parameter} is what distinguishes 
the real case from the complex case. We need some preparation to give this principle a precise meaning.

Let $R$ be a real closed field (for example, the field of real numbers or the field of algebraic real 
numbers).  A {\it real algebraic set} is an algebraic subset of $R^n$ for some $n$.
A {\it real algebraic manifold} is an irreducible nonsingular real algebraic set.
Given two real algebraic sets $X$ and $Y$, we write $X \sim Y$ to indicate that they are birationally 
isomorphic, $X\not\sim Y$ means that $X$ and $Y$ are not birationally isomorphic. 
In the following, we will use standard notions from semialgebraic and Nash geometry, 
as given in the book \lq\lq Real algebraic geometry\rq\rq\,  by Bochnak, Coste and Roy.

Let $f\colon X \rightarrow R^b$ be a surjective regular map defined on a real algebraic set $X$.
Set 
$$\mathcal S_f \ce \left\{ (y,z) \in R^b \times R^b: f^{-1}(y) = f^{-1}(z) \right\},$$
$$\rho_b\colon R^b \times R^b \rightarrow R^b, \quad (y,z) \mapsto  y.$$

The map $f$ is said to be {\it perfectly parametrized} by $R^b$ if
$f^{-1}(y) \not\sim f^{-1}(z) $ for all $y,z\in R^b$ with $y\neq z$; 
$f$ is said to be {\it almost perfectly parametrized} by $R^b$ if
there exists a semialgebraic subset $\mathcal T$ of $R^b\times R^b$
such that $\mathcal T$ contains $\mathcal S_f$  and each fiber of 
the restriction of $\rho_b$ to $\mathcal T$ is finite. 

Let $\pi\colon \mathcal V \rightarrow R^b$ be a surjective submersive regular map
from a real algebraic manifold $\mathcal V$ to $R^b$ with irreducible fibers. 
Such a map is called a {\it real algebraic family of real algebraic manifolds
parameterized by $R^b$}.

Given a real algebraic manifold $V$, we say that $\pi\colon \mathcal V \rightarrow R^b$ 
is an algebraic  real-deformation of $V$ 
if there exists a regular map $\pi' \colon \mathcal V \rightarrow V$ such that the map
$\pi \times \pi' \colon \mathcal V \rightarrow R^b \times V$ is a Nash isomorphism
and the restriction of $\pi'$ to $\pi^{-1}(0)$ is a biregular isomorphism.

The main result of \cite{13} is as follows.

\begin{theorem*}\cite{13} Every real algebraic manifold $V$ of positive dimension
has the following property: for each nonnegative integer $b$ there exists an algebraic real 
deformation of $V$ almost perfectly parameterized by $R^b$.
\end{theorem*}

In \cite{13} the authors also conjectured that in this theorem one can replace 
\lq\lq almost perfectly parameterized\rq\rq\,  by \lq\lq perfectly parameterized\rq\rq.
This conjecture was proved in a suitably modified form in \cite{14}. 

\section{Methodology and how to play the Ballico Game}

Today I must stop writing this text because  the deadline given by the editors has arrived. 
I would like to explain my methodology for arriving at this extremely short list of references.
When I  got the assignment of writing about the work of Ballico 
on moduli of vector bundles, my first search for articles authored by Ballico 
having  \lq\lq vector bundles\rq\rq\, resulted in  430 publications. Next, I refined the search, 
asking for 
\lq\lq vector bundles and moduli\rq\rq\, and it resulted in 324 publications by Ballico. 
 I was asked to start with some joint work, and this took me to applications to Physics, and
 upon reflexion, I decided to discuss those articles of Ballico on moduli spaces
 that give explicit contributions to mathematical physics, 
 plus two specific geometric
 features of moduli spaces, namely, the counting of irreducible components and 
 the existence of singularities. I would have liked to write more about the work of 
 Ballico on   \lq\lq deformations\rq\rq, but, even though this search resulted in \emph{only} 
 74 publications, it was still more than I could handle. So, I chose to present only 
 a result about moduli flexibility for real manifolds, because I found it very interesting. 
 However,  
 real varieties are outside of my area of expertise, so I cheated, asking Kucharz to give me the 
 text for this section. 
 
 So, this explains the choice of topics and also the bibliography included here.  The curious reader who 
 searches for the list of publications of Ballico having the word moduli of the title, can have the fun of verifying  
 that such list of articles is more than twice as long as this entire writeup, so that, just including those 
 references would already triplicate the length of this text. 
 
 Now that we are on the subject of number of publications, let me finish by telling something 
 about Edoardo that I think he does not know. For well over a decade we have been 
 having fun by playing
 the Ballico Game. \\
 
 \noindent{\bf The origin of the Ballico Game}: At some point,  more than 300 publications ago, 
 our group of about 10  mathematicians in Brasil had been looking at Ballico's papers very often. 
Since we were all very fond of Ballico's work, we wanted to congratulate him on the day 
when  his 1000th publication  appeared in MathSciNet. (Like Pel\'e was celebrated when he made his 1000th goal.)
 On Friday we checked and Ballico had 998 publications. We did not have access to MathSciNet
 during the weekend, 
 and on Monday morning we went to check it, Ballico had 1002 publications. We were both thrilled and disappointed, 
 because we  had missed the opportunity of sending the congratulations: \\
 
-- Goooool! Congratulations on the 1000th...\\

However, not all was lost. Since then we play the Ballico Game very often. In fact, just 
last week I played with 
Erica Klarreich, who I met  online on the occasion of an interview about Simons Associateships of ICTP. Erica won.  \\

 \noindent{\bf Rules of the Ballico Game}: The game may be played by any group of at least 2 players. 
 Each person picks a number. For the basic game, the referee for the game then  looks
  up the current number of Ballico's publications in MathSciNet 
 and the person whose guess is closest to the actual number wins. 
 
 For the advanced version of the game, useful to those wishing to play more than once in the same day
 (and we have done it),
 the referee picks a random date, and the goal is to  guess the  number of Ballico's publications prior to the
 chosen date.

\section{ Acknowledgements} 
E.G.  is a senior associate the Abdus Salam International Centre for Theoretical Physics, Italy, 
and thanks the ICTP for the hospitality during her visit in 2025, when the final version of this text was completed.

 \section{Conflict of Interest}
 The  author states that there is no conflict of interest.


\begin{thebibliography}{999}

\bibitem[1987]{87} Ballico, E.: {\it On moduli of vector bundles on rational surfaces},
Arch. Math. (Basel) {\bf 49} n.3  (1987)  267–272.%

\bibitem[1992]{92} Ballico, E.: {\it  On moduli of vector bundles on surfaces with negative Kodaira dimension},
Ann. Univ. Ferrara Sez. VII (N.S.) {\bf 38} (1992) 33–40.%

\bibitem[1993]{93} Ballico, E., Chiantini, L.: {\it On some moduli spaces of rank  2  bundles over  K3  surfaces},
Boll. Un. Mat. Ital. A {\bf 7} n.2 (1993) 279–287.%


\bibitem[1994]{94} Ballico, E.:
{\it On the components of moduli of simple sheaves on a surface},
Forum Math. {\bf 6} n.1 (1994) 35–41.%

\bibitem[1997a]{97a} Ballico, E.:
{\it On the number of components of the moduli schemes of stable torsion-free sheaves on integral curves},
Proc. Amer. Math. Soc. {\bf 125} n.10  (1997) 2819–2824.%

\bibitem[1997b]{97b} Ballico, E., Russo, B.: 
{\it  On the rationality of the moduli schemes of vector bundles on curves},
Rend. Sem. Mat. Univ. Padova {\bf 97} (1997) 23–27.%

\bibitem[1998a]{98a} Ballico, E.:
{\it  Singular moduli spaces of rank  2  stable vector bundles},
Quart. J. Math. Oxford Ser. (2) {\bf 49}  n.193 (1998) 31–34.%

\bibitem[1998b]{98b} Ballico, E.: 
{\it  On the birational structure of exceptional irreducible components of 
moduli schemes of rank  2  vector bundles on algebraic surfaces with negative Kodaira dimension},
Bull. Polish Acad. Sci. Math. {\bf 46} n.1 (1998) 17--23.%

\bibitem[1999]{99} Ballico, E.; Miró-Roig, R. M., 
{\it A lower bound for the number of components of the moduli schemes of stable rank  2  vector bundles on projective  3 -folds},
Proc. Amer. Math. Soc. {\bf 127} n.9 (1999) 2557--2560.%

\bibitem[2000]{00} Ballico, E.: {\it On the connectedness of the real part of moduli spaces of vector bundles on real algebraic surfaces},
J. Austral. Math. Soc. Ser. A {\bf 68} n.1 (2000) 41--54.

\bibitem[2003]{03} Ballico, E., Gasparim, E.:
{\it Vector  bundles on a neighborhood of an exceptional curve and elementary transformations},
Forum Math. {\bf 15} n.1 (2003) 115–122.%

\bibitem[2005]{05} Ballico, E., Gasparim, E,:
{\it Vector bundles on a three-dimensional neighborhood of a ruled surface},
J. Pure Appl. Algebra {\bf 195} n.1  (2005) 7–19.%

\bibitem[2009a]{09a} Ballico, E.; Gasparim, E.; Köppe, T.:
{\it Local moduli of holomorphic bundles},
J. Pure Appl. Algebra {\bf 213} n.4  (2009)  397--408.%

\bibitem[2009b]{09b} Ballico, E.; Gasparim, E.; Köppe, T.:
{\it Vector bundles near negative curves: moduli and local Euler characteristic},
Comm. Algebra {\bf 37} n.8  (2009) 2688--2713.%

\bibitem[2012]{12} Ballico, E., Eyral, C., Gasparim, E.:
{\it On the geometry of moduli spaces of anti-self-dual connections},
Topology Appl. {\bf 159} n.3 (2012) 633–645.%

\bibitem[2013]{13} Ballico, E., Ghiloni, R.:
{\it The principle of moduli flexibility for real algebraic manifolds},
Ann. Polon. Math. {\bf 109} n.1 (2013) 1--28.

\bibitem[2014]{14} Ballico, E., Ghiloni, R.:
{\it On the principle of real moduli flexibility: perfect parametrizations},
Ann. Polon. Math.   {\bf 111} n.3  (2014) 245--258.

\bibitem[2017]{17} Ballico, E., Gasparim, E., Grama, L., San Martin, L. A. B.:
{\it Some Landau-Ginzburg models viewed as rational maps},
Indag. Math. (N.S.) {\bf 28} n.3 (2017) 615--628.%


\bibitem[2018a]{18a} Ballico, E., Callander, B., Gasparim, E.:
{\it  Compactifications of adjoint orbits and their Hodge diamonds},
J. Algebra Appl. {\bf 17} n.6 (2018) 1850099, 16 pp.%

\bibitem[2018b]{18b} Ballico, E., Barmeier, S., Gasparim, E., Grama, L., San Martin, L.A.B.:
{\it A Lie theoretical construction of a Landau-Ginzburg model without projective mirrors},
Manuscripta Math. {\bf 158} n.1-2 (2019) 85–101.%

\bibitem[2021]{21} Ballico, E., Gasparim, E., Suzuki, B.: 
{\it Infinite dimensional families of Calabi--Yau threefolds and moduli of vector bundles},
J. Pure Appl. Algebra {\bf 225} n.4 (2021), 106554, 24 pp.%

\bibitem[2022a]{22a} Ballico, E., Gasparim, E., Rubilar, F., Suzuki, B.:
{\it The Kuranishi map for vector bundles on certain products of curves},
Topology Appl. {\bf 312} (2022), n. 108056, 12 pp.%

\bibitem[2022b]{22b} Ballico, E.; Gasparim, E.; Rubilar, F.:
{\it 25 open questions about vector bundles and their moduli},
Proyecciones {\bf 41} n.2 (2022) 403–436.%


\bibitem[2023a]{23a} Ballico, E., Gasparim, E., Köppe, T., Suzuki, B.:
Poisson structures on the conifold and local Calabi-Yau threefolds
Rep. Math. Phys. {\bf 90} n.3 (2022)  299–324.%

\bibitem[2023b]{23b} Ballico, E., Gasparim, E., Rubilar, F., Suzuki, B.:
{\it Lagrangian skeleta, collars and duality},
Springer Proc. Math. Stat. {\bf 409},
Springer, Cham, 2023, 263--280. %

\bibitem[2023c]{23c}  Ballico, E.,  Gasparim, E., Rubilar, F., Suzuki, B.:  
{\it Quantizations of local Calabi--Yau threefolds and their moduli of vector bundles},
 arXiv:2301.04192.%

\bibitem[2023d]{23d} Ballico, E., Gasparim, E.,  Rubilar, F.,  San Martin, L. A. B.:
{\it The Katzarkov--Kontsevich--Pantev conjecture for minimal adjoint orbits},
Eur. J. Math. {\bf 9} n.3  (2023) 57, 17 pp.%


















\end{thebibliography}
\end{document}